\documentclass[11pt]{article}   
\usepackage{amsthm, amsmath}   
\usepackage{amssymb}
\usepackage[hypertex]{hyperref}
\usepackage{algorithmic}
\usepackage{multirow}

\DeclareMathOperator{\Del}{Delta}

\theoremstyle{plain}

\newtheorem{theorem}{Theorem}[section]

\newtheorem{proposition}[theorem]{Proposition}
\newtheorem{corollary}[theorem]{Corollary}

\theoremstyle{definition}

\newtheorem{algorithm}[theorem]{Algorithm}

\theoremstyle{remark}

\title{The $3$-dimensional planar assignment problem and the
number of Latin squares related to an autotopism\footnote{Official printed version avalaible in Proceedings of XI Spanish Meeting on Computational Algebra and Applications EACA 2008 (2008), pp. 89-92.}.}
\author{
R. M. Falc\'on
\and
J. Mart\'in-Morales
}
\date{}

\begin{document}
\maketitle

\begin{abstract}
There exists a bijection between the set of Latin squares of order $n$ and the set of feasible solutions of the 3-dimensional planar assignment problem ($3PAP_n$). In this paper, we prove that, given a Latin square isotopism $\Theta$, we can add some linear constraints to the $3PAP_n$ in order to obtain a $1-1$ correspondence between the new set of feasible solutions and the set of Latin squares of order $n$ having $\Theta$ in their autotopism group. Moreover, we use Gr\"obner bases in order to describe an algorithm that allows one to obtain the cardinal of both sets.
\end{abstract}

\section*{Introduction}

A {\em Latin square} of order $n$ is an $n \times n$ array
with elements chosen from a set of $n$ distinct symbols (in this paper,
it will be the set $[n]=\{1,2,...,n\}$) such that each symbol occurs precisely once in
each row and each column. The set of Latin squares of order $n$ is
denoted by $LS(n)$. A {\em partial Latin square} of order
$n$, is a $n \times n$ array with elements chosen from a set of
$n$ symbols, such that each symbol occurs at most once in each row
and in each column. The set of partial Latin squares of order $n$
is denoted by $PLS(n)$.

The permutation group on $[n]$ is denoted by $S_n$. Every permutation $\delta\in S_n$ can be uniquely written as a
composition of $\mathbf{n}_{\delta}$ pairwise disjoint cycles, $\delta=C^{\delta}_1\circ C^{\delta}_2\circ ... \circ
C^{\delta}_{\mathbf{n}_{\delta}},$ where for all $i\in [\mathbf{n}_{\delta}]$, one has
$C^{\delta}_i=\left(c_{i,1}^{\delta}\ c_{i,2}^{\delta}\ ...\
c_{i,\ \lambda_i^{\delta}}^{\delta}\right)$, with
$c_{i,1}^{\delta}=\min_j \{c_{i,j}^{\delta}\}$. The {\em cycle structure of} $\delta$ is the sequence
$\mathbf{l}_{\delta}=(\mathbf{l}_1^{\delta},\mathbf{l}_2^{\delta},...,\mathbf{l}_n^{\delta})$,
where $\mathbf{l}_i^{\delta}$ is the number of cycles of length
$i$ in $\delta$, for all $i\in [n]$. Thus, $\mathbf{l}_1^{\delta}$ is the cardinal of the set of {\em fixed points} of $\delta$, $Fix(\delta)=\{i\in [n] \,\mid \, \delta(i)=i\}$. An {\em isotopism} of a Latin square $L\in LS(n)$ is a triple
$\Theta=(\alpha,\beta,\gamma)\in \mathcal{I}_n=S_n\times S_n\times
S_n$. Therefore, $\alpha,\beta$ and $\gamma$ are permutations of
rows, columns and symbols of $L$, respectively. The {\em cycle structure} of $\Theta$ is the triple
$(\mathbf{l}_{\alpha},\mathbf{l}_{\beta},\mathbf{l}_{\gamma})$.

An isotopism which maps $L$ to itself is an {\em autotopism}. The
possible cycle structures of the set of non-trivial autotopisms of
Latin squares of order up to $11$ were obtained in
\cite{FalconAC}. The set of all possible autotopisms of order $n$ is denoted by $\mathfrak{A}_n$.
The stabilizer subgroup of $L$ in $\mathfrak{A}_n$
is its {\em autotopism group} $\mathfrak{A}(L)$. Given
$\Theta\in\mathfrak{A}_n$, the set of all Latin squares $L$ such
that $\Theta\in \mathfrak{A}(L)$ is denoted by $LS(\Theta)$ and
the cardinality of $LS(\Theta)$ is denoted by $\Delta(\Theta)$.
Specifically, if $\Theta_1$ and $\Theta_2$ are two autotopisms
with the same cycle structure, then
$\Delta(\Theta_1)=\Delta(\Theta_2)$. Given $\Theta\in
\mathfrak{A}_n$ and $P\in PLS(n)$, the number
$c_P=\Delta(\Theta)/|LS_P(\Theta)|$ is called the {\em $P$-coefficient
of symmetry of} $\Theta$, where $LS_P(\Theta)=\{L\in LS(\Theta)
\,\mid \, P\subseteq L\}$.

Gr\"obner bases were used in \cite{FalconMartinJSC07} to describe
an algorithm that allows one to obtain the number $\Delta(\Theta)$
in a computational way. This algorithm was implemented in {\sc
Singular} \cite{Greuel05} for Latin squares of
order up to $7$ \cite{FalconMartinLS07}. However, after applying it to upper orders, the authors have seen that,
in order to improve the time of computation, it is convenient to
combine Gr\"obner bases with some combinatorial tools. In this paper we study, as a possible tool, the 1-1 correspondence between $LS(n)$ and the set of feasible solutions of the {\em 3-dimensional planar assignment problem} ($3PAP_n$) \cite{Euler86}:
{\small $$\min \sum_{i\in I, j\in J, k\in K} w_{ijk}\cdot x_{ijk}, \hspace{1cm} s.t. \begin{cases}
         \sum_{i\in I}x_{ijk}=1, \forall j\in J, k\in K. \\
         \sum_{j\in J}x_{ijk}=1, \forall i\in I, k\in K. \\
         \sum_{k\in K}x_{ijk}=1, \forall i\in I, j\in J. \\
         x_{ijk}\in\{0,1\}, \forall i\in I, j\in J, k \in K.
       \end{cases} \hspace{2cm} (1)
$$}
where $w_{ijk}$ are real weights and $I,J,K$ are three disjoint $n$-sets. Thus, any feasible solution of the $3PAP_n$ can be considered as a Latin square $L=(l_{i,j})\in LS(n)$, by taking $I=J=K=[n]$ and $x_{ijk}=\begin{cases} 1, \text{ if } l_{i,j}=k,\\ 0, \text{ otherwise}.\end{cases}$. The reciprocal is analogous.

\section{Constraints related to an autotopism of a Latin square}

Given a autotopism $\Theta=(\alpha,\beta,\gamma)\in \mathfrak{A}_n$, let $(1)_{\Theta}$ be the set of constraints obtained by adding to $(1)$ the $n^3$ constraints $x_{ijk}=x_{\alpha(i)\beta(j)\gamma(k)}$, $\forall i\in I, j\in J, k\in K$.

\begin{theorem} \label{thr0}
There exists a bijection between $LS(\Theta)$ and the set of feasible solutions related to a combinatorial optimization problem having $(1)_{\Theta}$ as the set of constraints.  \hfill $\Box$
\end{theorem}

$(1)_{\Theta}$ is a system of $3n^2 + 2n^3$ equations of degrees $1$ and $2$, in $n^3$ variables, which can be solved by using Gr\"obner basis. Thus, if we define $F(x)=x\cdot(x-1)$, then the following result is verified:

\begin{corollary} \label{crl1} $LS(\Theta)$ corresponds to the set of zeros of the ideal $I = \langle\, \sum_{i\in [n]}x_{ijk}-1 \mid j,k
 \in [n]\rangle\, +
      \langle\, \sum_{j\in [n]}x_{ijk}-1 \mid i,k \in [n]\rangle\, +
      \langle\, \sum_{k\in [n]}x_{ijk}-1 \mid i,j \in [n]\rangle\, + \langle\, F(x_{ijk}) \mid i,j,k\in [n]\rangle\, +
      \langle\, x_{ijk}-x_{\alpha(i)\beta(j)\gamma(k)} \mid i,j,k \in [n]\rangle\,\subseteq \mathbb{Q}[{\bf x}] = \mathbb{Q}[x_{111},...,$ $x_{nnn}]$. $\hfill \Box$
\end{corollary}

The symmetrical structure of $\Theta$ can be used to reduce the number of variables of the previous system. To see it, let us consider {\small
$S_{\Theta}  = \left\{(i,j) \,\mid \, i\in S_{\alpha},j\in \begin{cases} [n],
     \text{ if }i\not \in Fix(\alpha),\\
     S_{\beta}, \text{ if } i\in Fix(\alpha).\end{cases}\right\}$} as a set of $(\mathbf{n}_{\alpha}-\mathbf{l}_{\alpha}^1)\cdot n + \mathbf{l}_{\alpha}^1\cdot \mathbf{n}_{\beta}$ multi-indices, where $S_{\alpha}=\{c^{\alpha}_{i,1}\mid i\in [\mathbf{n}_{\alpha}]\}$ and $S_{\beta}=\{c^{\beta}_{j,1}\mid j\in [\mathbf{n}_{\beta}]\}$.

\begin{proposition}[Falc\'on and Mart\'\i n-Morales \cite{FalconMartinJSC07}] \label{prp1} Let $L=(l_{i,j})\in LS(\Theta)$ be such that
all the triples of the Latin subrectangle $R_L=\left\{(i,j,l_{i,j})\,\mid (i,j)\in S_{\Theta}\right\}$ of $L$ are known. Then, all the triples of $L$ are known.\hfill $\Box$
\end{proposition}

Let $\varphi_{\Theta}$ be a map in the set of $n^3$ variables ${\bf x}=\{x_{111},...,x_{nnn}\} $ such that $\varphi_{\Theta}(x_{ijk})=\begin{cases} x_{ijk}, \text{ if } (i,j)\in S_{\Theta},\\
x_{\alpha^m(i)\beta^m(j)\gamma^m(k)}, \text{ otherwise}.
\end{cases},$ where $m=\min \{l\in [n] \mid
(\alpha^l(i),\beta^l(j))\in S_{\Theta}\}$.

\begin{theorem} \label{thr1} $LS(\Theta)$ corresponds to the set of zeros of the ideal $I'= \langle\, \sum_{i\in [n]}\varphi_{\Theta}(x_{ijk})-1 \mid j,k \in [n]\rangle\, +
      \langle\, \sum_{j\in [n]}\varphi_{\Theta}(x_{ijk})-1 \mid i,k \in [n]\rangle\, +
      \langle\, \sum_{k\in [n]}\varphi_{\Theta}(x_{ijk})-1 \mid i,j \in [n]\rangle\, +
      \langle\, x_{ijk} \mid \alpha(i)=i, \beta(j)=j, \gamma(k)\neq k \rangle\, +
      \langle\, F(x_{ijk}) \mid (i,j)\in S_{\Theta}, k\in [n] \rangle\,=
      \langle\, \varphi_{\Theta}(I)\rangle\,\subseteq \mathbb{Q}[\varphi_{\Theta}({\bf x})]$. $\hfill \Box$
\end{theorem}

  Now, let $P=(p_{i,j})\in PLS(n)$ be such that
  $p_{i,j}=\emptyset$, for all $(i,j)\not\in S_{\Theta}$ and let $c_P$ be the $P$-coefficient
  of symmetry of $\Theta$. Thus, we know that $\Delta(\Theta)=c_P\cdot |LS_P(\Theta)|$
  and we can calculate $|LS_P(\Theta)|$ starting from the set of solutions
  of an algebraic system of polynomial equations associated with $\Theta$ and $P$. Specifically,
  we obtain the following algorithm:

  \begin{algorithm}[Computation of $\Delta(\Theta)$]
  \label{alg1}
  \begin{algorithmic}\ \\
    \STATE Input: \ $\Theta = (\alpha,\beta,\gamma)\in\mathcal{I}_n$;\\
    \hspace{1.22cm} $\mathbf{n}_\alpha$, the number of cycles of $\alpha$;\\
    \hspace{1.22cm} $P\in PLS(n)$ such that
    $p_{i,j}=\emptyset$, for all $(i,j)\not\in S_{\Theta}$;\\
    \hspace{1.22cm} $c_P$, the $P$-coefficient of symmetry of $\Theta$.\smallskip
    \STATE Output: \ $\Delta(\Theta)$, the number of Latin squares having $\Theta$ as an autotopism;\bigskip
    \STATE {\small $I':= \langle\, \sum_{i\in [n]}\varphi_{\Theta}(x_{ijk})-1 \mid j,k \in [n]\rangle\, +
      \langle\, \sum_{j\in [n]}\varphi_{\Theta}(x_{ijk})-1 \mid i,k \in [n]\rangle\, +
      \langle\, \sum_{k\in [n]}\varphi_{\Theta}(x_{ijk})-1 \mid i,j \in [n]\rangle\, + \langle\, F(x_{ijk}) \mid (i,j)\in S_{\Theta}, k\in [n] \rangle$};\\[0.2cm]
    \STATE $I':=I' \, + \,\langle\, x_{ijl}-\delta^l_{p_{i,j}} \mid p_{i,j}\neq \emptyset, l\in [n]\,\rangle\, + \,\langle\, x_{ilp_{i,j}}-\delta^l_j \mid  p_{i,j}\neq \emptyset, l\in [n]\,\rangle\, + \,\langle\, x_{ljp_{i,j}}-\delta^l_i \mid  p_{i,j}\neq \emptyset, l\in [n]\,\rangle\,$;
    \hfill $\triangleright$ $\delta$ is Kronecker's delta.\\
    \STATE $G:=$ Gr\"obner basis of $I'$ with respect to any term ordering;
    \STATE $\Del := \dim_\mathbb{Q} (\mathbb{Q}[\varphi_{\Theta}({\bf x})]/I')$;
    \hfill $\triangleright$ $\Del$ is the cardinality of $V(I')$\\
    RETURN $c_P\cdot\Del$;
  \end{algorithmic}
  \end{algorithm}

\section{Number of Latin squares related to $\mathfrak{A}_8$ and
$\mathfrak{A}_9$.}

We have implemented Algorithm \ref{alg1} in a {\sc
Singular} procedure \cite{FalconMartinEACA08} which
improves running times of \cite{FalconMartinJSC07}. Moreover, we have obtained the number
$\Delta(\Theta)$ corresponding to autotopisms of
$\mathfrak{A}_8$ and $\mathfrak{A}_9$, as we can see in Table 1. The timing information, measured in
seconds, has been taken from an {\em Intel Core 2 Duo Processor T5500, 1.66 GHz}
with {\em
Windows Vista} operating system.

\begin{table}[ht]
\centering{\tiny
\begin{tabular}{|c|c|c|c|c|}\hline
$n$ & $\mathbf{l}_{\alpha} = \mathbf{l}_{\beta}$ &
$\mathbf{l}_{\gamma}$ & $\Delta$ & \begin{tabular}{c} r.t. \\
                                                                    (\ref{alg1})
                                                                    \end{tabular}\\
\hline \multirow{5}{*}{8} \ & \ & (0,0,0,2,0,0,0,0) & 1152 & 16\\
\cline{3-5}
\ & \ & (0,2,0,1,0,0,0,0) & 1408 & 12\\
\cline{3-5}
\ & \ & (0,4,0,0,0,0,0,0) & 3456 & 10\\
\cline{3-5}
\ & \ & (2,1,0,1,0,0,0,0) & 1408 & 14\\
\cline{3-5}
\ & (0,0,0,0,0,0,0,1) & (2,3,0,0,0,0,0,0) & 3456 & 13\\
\cline{3-5}
\ & \ & (4,0,0,1,0,0,0,0) & 3456 & 15\\
\cline{3-5}
\ & \ & (4,2,0,0,0,0,0,0) & 8064 & 21\\
\cline{3-5}
\ & \ & (6,1,0,0,0,0,0,0) & 17280 & 34\\
\cline{3-5}
\ & \ & (8,0,0,0,0,0,0,0) & 40320 & 12\\
\cline{2-5}
\ & \ & (0,0,0,2,0,0,0,0) & 106496 & 945\\
\cline{3-5}
\ & \ & (0,2,0,1,0,0,0,0) & 188416 & 1163 \\
\cline{3-5}
\ & \ & (0,4,0,0,0,0,0,0) & 811008 & 255 \\
\cline{3-5}
\ & \ & (2,1,0,1,0,0,0,0) & 253952 & 731 \\
\cline{3-5}
\ & (0,0,0,2,0,0,0,0) & (2,3,0,0,0,0,0,0) & 1007616 & 548 \\
\cline{3-5}
\ & \ & (4,0,0,1,0,0,0,0) & 712704 & 600 \\
\cline{3-5}
\ & \ & (4,2,0,0,0,0,0,0) & 2727936 & 660 \\
\cline{3-5}
\ & \ & (6,1,0,0,0,0,0,0) & 7741440 & 73\\
\cline{3-5}
\ & \ & (8,0,0,0,0,0,0,0) & 23224320 & 1\\
\cline{2-5}
\ & (0,1,0,0,0,1,0,0) & (2,0,0,0,0,1,0,0) & 3456 & 5 \\
\cline{3-5}
\ & \ & (2,0,2,0,0,0,0,0) & 19008 & 3 \\
\cline{2-5}
\ & (1,0,0,0,0,0,1,0) & (1,0,0,0,0,0,1,0) & 931 & 76 \\
\cline{2-5}
\ & \ & (0,2,0,1,0,0,0,0) & 16384 & 3 \\
\cline{3-5}
\ & (0,2,0,1,0,0,0,0) & (2,1,0,1,0,0,0,0) & 16384 & 3 \\
\cline{3-5}
\ & \ & (4,0,0,1,0,0,0,0) & 147456 & 3 \\
\cline{2-5}
\ & (2,0,0,0,0,1,0,0) & (2,0,0,0,0,1,0,0) & 19584 & 72 \\
\cline{3-5}
\ & (0,4,0,0,0,0,0,0) & (6,1,0,0,0,0,0,0) & 198747095040 & 6515\\
\cline{3-5}
\ & \ & (8,0,0,0,0,0,0,0) & 828396011520 & 9027\\
\cline{2-5}
\ & (2,1,0,1,0,0,0,0) & (2,1,0,1,0,0,0,0) & 8192 & 1\\
\cline{2-5}
\ & (3,0,0,0,1,0,0,0) & (3,0,0,0,1,0,0,0) & 388800 & 80\\
\cline{2-5}
\ & (4,0,0,1,0,0,0,0) & (4,0,0,1,0,0,0,0) & 7962624 & 2\\
\cline{2-5}
\ & (4,2,0,0,0,0,0,0) & (4,2,0,0,0,0,0,0) & 509607936 & 10\\
\hline \multirow{5}{*}{9} \ & \ & (0,0,0,0,0,0,0,0,1) & 2025 & 50\\
\cline{3-5} \ & \ & (0,0,3,0,0,0,0,0,0) & 7128 & 33\\
\cline{3-5} \ & (0,0,0,0,0,0,0,0,1) & (3,0,2,0,0,0,0,0,0) & 12960 & 61\\
\cline{3-5} \ & \ & (6,0,1,0,0,0,0,0,0) & 71280 & 221\\
\cline{3-5} \ & \ & (9,0,0,0,0,0,0,0,0) & 362880 & 3\\
\cline{2-5} \ & \ & (0,0,1,0,0,1,0,0,0) & 15552 & 46\\
\cline{3-5} \ & (0,0,1,0,0,1,0,0,0) & (0,3,1,0,0,0,0,0,0) & 124416 & 4\\
\cline{3-5} \ & \ & (3,0,0,0,0,1,0,0,0) & 62208 & 16\\
\cline{3-5} \ & \ & (3,3,0,0,0,0,0,0,0) & 1244160 & 17\\
\cline{2-5} \ & (1,0,0,0,0,0,0,1,0) & (1,0,0,0,0,0,0,1,0) & 4096 & 56\\
\cline{2-5} \ & (0,0,3,0,0,0,0,0,0) & (6,0,1,0,0,0,0,0,0) & 403813278720 & 221\\
\cline{3-5} \ & \ & (9,0,0,0,0,0,0,0,0) & 948109639680 & 1846\\
\cline{2-5} \ & (1,0,0,2,0,0,0,0,0) & (1,0,0,2,0,0,0,0,0) & 12189696 & 11098\\
\cline{2-5} \ & (1,1,0,0,0,1,0,0,0) & (1,1,0,0,0,1,0,0,0) & 69120 & 557\\
\cline{2-5} \ & (2,0,0,0,0,0,1,0,0) & (2,0,0,0,0,0,1,0,0) & 438256 & 615\\
\cline{2-5} \ & (3,0,0,0,0,1,0,0,0) & (3,0,0,0,0,1,0,0,0) & 3110400 & 112\\
\cline{2-5} \ & (4,0,0,0,1,0,0,0,0) & (4,0,0,0,1,0,0,0,0) & 199065600 & 3\\
\hline
\end{tabular}}\bigskip\\
\caption{Number of Latin squares related to $\mathfrak{A}_8$ and
$\mathfrak{A}_9$.}
\end{table}

{\small

\begin{tabular}{lccl}
{\em R. M. Falc\'on} & \ & \ & {\em J. Mart\'in-Morales}\\
Department of Applied Mathematics I. & \ & \ & Department of Mathematics.\\
Technical Architecture School. University of Seville. & \ & \ & University of Zaragoza.\\
Avda. Reina Mercedes, 4A. & \ & \ & C/ Pedro Cerbuna, 12.\\
41012, Seville (Spain). & \ & \ & 50009, Zaragoza (Spain).\\
E-mail: {\em rafalgan@us.es}& \ & \ & E-mail: {\em jorge@unizar.es}
\end{tabular}}


\begin{thebibliography}{}
{\small

\bibitem{Adams94} Adams, W. and
Loustaunau, P., 1994. \newblock \emph{An introduction to Gr\"obner
bases}. \newblock Volume 3 of Graduate Studies in Mathematics.
American Mathematical Society, Providence, RI.

\bibitem{Euler86} Euler, R., Burkard, R. E. and Grommes, R. \newblock \emph{On Latin squares and the facial structure of related polytopes}.
\newblock Discrete Mathematics 62 (1986), pp. 155 - 181.

\bibitem{FalconAC} Falc\'on, R. M. \newblock \emph{Cycle structures of autotopisms of the Latin squares of order
up to 11}. \newblock Ars Combinatoria (in press). \newblock Avalaible from
\url{http://arxiv.org/abs/0709.2973}.

\bibitem{FalconMartinJSC07} Falc\'on, R. M. and Mart\'in-Morales, J.
\newblock \emph{Gr\"obner bases and the number of Latin squares related to autotopisms of
order $\leq$ 7}. \newblock Journal of Symbolic Computation 42
(2007), pp. 1142 - 1154.

\bibitem{FalconMartinLS07}
\url{http://www.personal.us.es/raufalgan/LS/latinSquare.lib}

\bibitem{FalconMartinEACA08}
\url{http://www.personal.us.es/raufalgan/LS/3PAPlatinSquare.lib}

\bibitem{Greuel05}
Greuel, G.-M., Pfister, G. and Sch\"onemann, H., 2005. \newblock
{\sc Singular} 3.0. A Computer Algebra System for Polynomial
Computations. \newblock Centre for Computer Algebra, University of
Kaiserlautern. \newblock \url{http://www.singular.uni-kl.de}.}

\end{thebibliography}
\end{document}